\documentclass[11pt,twoside]{article}
\usepackage[dvips]{graphicx}
\usepackage{fancyhdr}
\usepackage{amssymb, amsmath,color,graphicx, subfigure, cite, tikz, sidecap}

\usetikzlibrary{arrows,calc}
\tikzset{
>=stealth',
help lines/.style={dashed, thick},
axis/.style={<->},
important line/.style={thick},
connection/.style={thick, dotted},
}
\newtheorem{thm}{Theorem}[section]

\newtheorem{lem}[thm]{Lemma}
\newtheorem{prop}[thm]{Proposition}

\newcommand{\norm}[1]{\left\Vert#1\right\Vert}
\newcommand{\abs}[1]{\left\vert#1\right\vert}
\newcommand{\set}[1]{\left\{#1\right\}}
\newcommand{\pa}[1]{\left(#1\right)}
\newcommand{\ma}[1]{\left[#1\right]}

\newcommand{\inne}[1]{\left\langle #1\right\rangle}

\newcommand{\reals}{{\mathbb{R}}}

\newcommand{\complexes}{{\mathbb{C}}}

\newcommand{\PROOF}{\noindent {\bf Proof}: }

\newcommand{\cT}{\ma{\begin{array}{rrrrr}
0  & 0 & ... & 1\\
1  & 0 & ... & 0\\
   & \ddots & \ddots & \\
0  & ...& 1 & 0

\end{array}}}

\newcommand{\EE}{\begin{bmatrix} \vline& \vline&&\vline\\
e_0 & e_1 & ... & e_{N-1}\\
\vline& \vline&&\vline\end{bmatrix}}

\usetikzlibrary[topaths]
\newcount\mycount

\begin{document}

%
%
\title{\bf\vspace{-39pt} Graph theoretic uncertainty and feasibility}

%
%

\author{Paul J. Koprowski\\  \small pkoprows@math.umd.edu} 


%
%

%
%
\maketitle
\thispagestyle{fancy}

%
%
\markboth{\footnotesize \rm \hfill PAUL J. KOPROWSKI \hfill}
{\footnotesize \rm \hfill Graph theoretic uncertainty and feasibility \hfill}

%
%

\begin{abstract} We expand upon a graph theoretic set of uncertainty principles with tight 
bounds for difference estimators acting simultaneously in the
graph domain and the frequency domain. We show that the eigenfunctions 
of a modified graph Laplacian and a modified normalized graph Laplacian
operator dictate the upper and lower bounds for the inequalities.
Finally, we establish the feasibility region of difference estimator values in $\reals^2$.
\vspace{5mm} 
\vspace{3mm}
\end{abstract}

\section{Introduction} Analysis on graphs is a key component to many 
techniques in data analysis, dimension reduction, and
analysis on fractals. The Fourier transform on a graph
has been defined using the spectrum of the graph Laplacian,
see, e.g., \cite{graphwavelets}, \cite{chungbook}, \cite{shuman2015vertex},
\cite{shuman2012windowed}, \cite{shuman2013emerging},
\cite{sandryhaila2014discrete}, \cite{pham2014wavelets}, \cite{ekambaram2013wavelet},
and \cite{agaskar2013spectral}.
In \cite{agaskar2013spectral}, the authors define the notion of spread 
in the spectral and graph domains using the analytic properties of the graph
Fourier transform. More recently in \cite{tsitsvero2015signals} and \cite{tsitsvero2015uncertainty},
the authors introduce the notion of band limiting operators and
the effects of such operations in the graph setting.
The eigenvalues and eigenvectors of the graph Laplacian 
play a central role in the uncertainty analysis in the aforementioned
papers.
Motivated by the feasibility results in \cite{agaskar2013spectral},
we extend the notion of discrete uncertainty 
principles such as those introduced in \cite{grunbaum2003heisenberg},
\cite{slepian1978prolate}, and \cite{donoho1989uncertainty}. We show 
that for the graph setting, the cyclic structure of the discrete Fourier
transform is no longer present for the graph Fourier transform. As a result, 
the support theorems (such as in \cite{donoho1989uncertainty}) are no longer
guaranteed. We extend the frame uncertainty principle introduced by 
Lammers and Maeser in \cite{lammers2011uncertainty}. Finally, we establish
the feasible values of the difference estimators, and provide explicit analysis
in the case of complete graphs. The bulk of this paper is an expansion of
the work published in \cite{benedetto2015graph} and detailed in \cite{koprowski2015finite}.

The structure of the paper is as follows. In Section \ref{graphs}, we provide
an overview of elementary graph theory, and we establish notation; 
and in Section \ref{graphFT} we define the graph Fourier transform.
Additive graph uncertainty principles are established in Section \ref{GUP}.
In Section \ref{GFUP}, we extend a result from \cite{lammers2011uncertainty} to the graph setting. 
Section \ref{feasibility} details the feasibility region in $\reals^2$ of possible
difference estimator values. 
In Section \ref{complete}, we provide uncertainty analysis on complete graphs.
Theorems \ref{Koprowski1}, \ref{Koprowski2}, \ref{graph_frame_up}, \ref{graph_frame_up_norm}, \ref{DNFR}
and \ref{DUC}
are the main results of the paper. The proof of Theorems \ref{DNFR} and
\ref{DUC}
follow the arguments by Agaskar and Lu in \cite{agaskar2013spectral}. As such, 
Section \ref{feasibility} is quite protracted.

\section{Weighted Graphs} \label{graphs}

A graph $G=\set{V,\mathbf{E}\subseteq V\times V,w}$ consists of a set $V$ of vertices, 
a set $\mathbf{E}$ of edges consisting of pairs of elements of $V,$ and a weight 
function $w: V\times V\to \reals^+$. For $u,v\in V$, $w(u,v)>0$ if $(u,v)\in \mathbf{E}$
and is zero otherwise. If $w(u,v)=1$ for all $(u,v)\in \mathbf{E}$, then we say $G$ is 
\textit{unit weighted} (or unweighted).
There is no restriction on the size of the set $V$, but we shall restrict our attention 
to $\abs{V}=N<\infty$. We also assume that the set $\set{v_j}_{j=0}^{N-1}=V$ has an arbitrary, but fixed ordering.

For all graphs, we define the \textit{adjacency} matrix $A$ component-wise as 
$A_{m,n}=w(v_m,v_n)$. If $A$ is symmetric, that is, if
$w(v_n,v_m)=A_{n,m}=A_{m,n}=w(v_m,v_n),$ then we say $G$ is \textit{undirected}.
If a graph has \textit{loops}, that is $w(v_j,v_j)>0$
for some $v_j\in V$, then $A$ has nonzero diagonal entries. 
Unless otherwise specified, we shall assume that our graphs are undirected, and have no loops.
We shall refer to such graphs as \textit{simple}.


The \textit{degree} $d$ of a vertex $v_j$ is defined by $deg(v_j)=\sum_{n=0}^{N-1}w(v_j,v_n)=\sum_{n=0}^{N-1}A_{j,n}$.
We can then define a diagonal \textit{degree matrix}, $$D=\mbox{diag}(deg(v_0), deg(v_1),...,deg(v_{N-1})).$$
There are two common choices for the graph Laplacian:
\begin{align*}\begin{array}{rcl}L&=&D-A\\
\mathcal{L}&=&I-D^{-1/2}AD^{-1/2},
\end{array}\end{align*}
where $I$ is the $N\times N$ identity. $L$ is the 
\textit{graph Laplacian}, while $\mathcal{L}$ is the
\textit{normalized graph Laplacian}.
Define the $|\mathbf{E}|\times |V|$ \textit{incidence matrix}
$M$ with element $M_{k,j}$ for edge $e_k$ and vertex $v_j$ by:
$$M_{k,j} = 
\begin{cases}\phantom{-} 1, & \text{if }\,e_k=(v_j,v_l)\mbox{ and } j<l\\
-1, & \text{if }\,e_k=(v_j,v_l)\mbox{ and } j>l\\
\phantom{-} 0, & \mbox{otherwise}.\end{cases}$$
Define the diagonal $\abs{\mathbf{E}}\times\abs{\mathbf{E}}$ 
\textit{weight matrix} $W=diag(w(e_0),w(e_1),...,w(e_{\abs{\mathbf{E}}-1})).$

Noting that $L=M^*WM=\pa{W^{\frac{1}{2}}M}^*\pa{W^{\frac{1}{2}}M}$,  where $\cdot^*$ denotes the conjugate
transpose of an operator $\cdot$, we conclude that $L$ is real, symmetric, and
positive semidefinite. By the spectral theorem, $L$ must have an 
orthonormal basis $\set{\chi_l}$ of eigenvectors with associated
eigenvalues $\set{\lambda_l}$ ordered as $0=\lambda_0<\lambda_1
\leq \lambda_2\leq...\leq\lambda_{N-1}.$ Let $\chi$ be the matrix whose $l^
{th}$ column is defined by $\chi_l$. Let $\Delta$ be the diagonalization 
of $L$, that is, $\chi^* L \chi=\Delta=diag(\lambda_0,...,
\lambda_{N-1}).$
We shall use this set of eigenfunctions $\set{\chi_l}$
to define the graph Fourier transform
in Section \ref{graphFT}.

Alternatively, after noting that $$\mathcal{L}=D^{-1/2}LD^{-1/2}=\pa{W^{\frac{1}{2}}MD^{-1/2}}^*\pa{W^{\frac{1}{2}}MD^{-1/2}},$$
we may apply the 
spectral theorem to $\mathcal{L}$. Hence, $\mathcal{L}$ must have an orthonormal eigenbasis $\set{F_l}$ with associated eigenvalues $\set{\mu_l}$ ordered as
$0=\mu_0<\mu_1\leq \mu_2\leq...\leq 
\mu_{N-1}.$ Let $\mathcal{F}$ be the matrix whose $l^{th}$ column is defined
by $F_l$ such that 
$\mathcal{F}$ diagonalizes $\mathcal{L}$. We shall use
this set of eigenfunctions $\set{F_l}$ to define the normalized graph Fourier transform
in Section \ref{graphFT}.



\section{The Graph Fourier Transform} \label{graphFT}
Functions $f$ defined on a graph $G$ will be written 
notationally as a vector $f\in\reals^N$ where
$f[j]$ for $j=0,...,N-1$ is the value of the function
$f$ evaluated at the vertex $v_j$.
We say $f\in l^2(G)$, and use the standard $l^2$
norm: $\norm{f}=\pa{\sum_{j=0}^{N-1}\abs{f[j]}^2}^{1/2}$.

Given this space $l^2(G)$ of real-valued functions on the set $V$ of vertices of the graph $G$,
it is natural to define a Fourier transform based on the structure of $G$. To motivate this definition, we examine the \textit{Fourier transform} on $L^1(\reals),$ viz.,
\begin{align*} \widehat{f}(\gamma)=\int_\reals f(t)e^{-2\pi i t\gamma}d\gamma, \end{align*}
and the formal \textit{inverse Fourier transform},
 \begin{align*}f(t)=\int_{\widehat{\reals}} \widehat{f}(\gamma)e^{2\pi it\gamma }\,d\gamma,\end{align*}
 where $\widehat{\reals}=\reals$ is considered the frequency domain. The functions $e^
{2\pi it\gamma }$, $\gamma\in\widehat{\reals},$ are the eigenfunctions of the Laplacian operator $\frac{d^2}{dt^2}$ since we have $\frac{d^2}{dt^2}e^{2\pi i t \gamma}
=-4\pi^2\gamma^2e^{2\pi i t \gamma}$. If $\widehat{f}\in L^1(\widehat{\reals})$, then the inverse Fourier transform is an expansion of the
function $f$ in terms of the eigenfunctions with coefficients $\widehat{f}(\gamma)$. With this in mind, we use the eigenvectors of the graph
Laplacian to
define the \textit{graph Fourier transform} $\widehat{f}$ of $f\in l^2(G)$ as follows:
$$\forall l=0,1,...,N-1,\quad \widehat{f}[l]=\inne{f,\chi_l},$$ or, equivalently, $\widehat{f}=\chi^* f$.
It is clear from the orthonormality of the basis, $\set{\chi_l}$, that $\chi^*=\chi^{-1}.$ Thus, the \textit{inverse graph Fourier transform} is given by $\chi \widehat
 {f}=\chi\chi^*f=If=f$, or, equivalently, $f[j]=\sum_{l=0}^{N-1}\inne{f,\chi_l}\chi_l[j]$.

Similarly, we define the \textit{normalized graph Fourier transform} $\overset{*}{f}$ of $f\in l^2(G)$ as follows:
$$\forall l=0,1,...,N-1,\quad \overset{*}{f}[l]=\inne{f,F_l},$$ or, equivalently, $\overset{*}{f}=\mathcal{F}^* f$.
It is clear from the orthonormality of the basis, $\set{F_l}$, that $\mathcal{F}^*=\mathcal{F}^{-1}.$ Thus, the \textit{inverse normalized graph Fourier 
transform} is given by $$\mathcal{F} \overset{*}{f}
=\mathcal{F}\mathcal{F}^*f=If=f,$$ or, equivalently, $f[j]=\sum_{l=0}^{N-1}\inne{f,F_l}F_l[j]$.

\section{Graph Uncertainty Principles}\label{GUP}
\begin{figure}[h]
  \centering
\begin{tikzpicture}

\def \n {7}
\def \radius {1.7cm}
\def \margin {8} 

\foreach \s in {0,...,\n}
{
  \node[draw, circle] at ({360/(\n+1) * (\s )}:\radius) {$\s$};
  \draw[-, >=latex] ({360/(\n+1) * (\s - 1)+\margin}:\radius)
    arc ({360/(\n+1) * (\s - 1)+\margin}:{360/(\n+1) * (\s)-\margin}:\radius);
}
\end{tikzpicture}
 \caption{\textbf{A unit weighted circulant graph with 8 vertices.}  The graph Laplacian associated with
 this graph is the classical discrete Laplacian.\vspace{1cm}}\label{circ}
\end{figure}
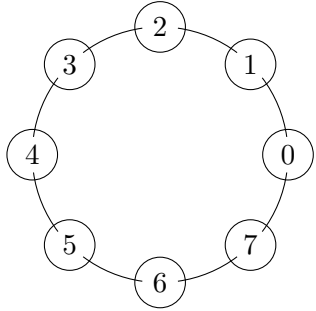
In the classical $L^2(\reals)$ setting, we have the additive Heisenberg uncertainty principle:
\begin{align} \norm{f(t)}^2\leq 2\pi\pa{\norm{tf(t)}^2+\norm{\gamma 
\widehat{f}(\gamma)}^2}\label{CUP}.\end{align}
For a function $f\in\mathcal{S}(\reals),$ the space of Schwartz functions on $\reals$, inequality (\ref{CUP}) is equivalent to:
\begin{align} \norm{f(t)}^2\leq \pa{\norm{ \widehat{f}'(\gamma)}^2+\norm{f'(t)}^2}.\label{CUP2}\end{align} To achieve a graph analog of inequality
(\ref{CUP2}), we must define the notion of a derivative or
difference operator in the graph setting. To do this, we examine the following product:
 $$W^{1/2}Mf=D_rf,$$ where $D_r=W^{1/2}M$. 
We refer to $D_r$ as the \textit{difference operator}
for the graph $G$ because it generates the weighted difference
of $f$ across each edge of $G$:
$$(D_rf)[k]=\pa{f[j]-f[i]}(w(e_k))^{1/2},$$
where $e_k=(v_j,v_i)$ and $j<i$.
It should be noted that the convention $j<i$ is arbitrary, and the
importantance of $D_r$ is in its magnitude. Indeed, 
in \cite{shuman2013emerging}, $\norm{D_r f}$ is the smoothness
measure of a function.
It is also common to refer to $D_rf$
as the derivative of $f$ (see \cite{agaskar2013spectral}).
In the case of the unit weighted circulant graph (see Figure \ref{circ}), $D_r$ is the
difference operator in \cite{lammers2011uncertainty}.
With this in mind, we establish a graph Fourier transform inequality of the form of (\ref{CUP2}).

\begin{thm}\label{Koprowski1}
Let $G$ be a simple, connected, and undirected graph. Then, for any non-zero function $f\in l^2 (G)$, the following inequalities hold:
\begin{align} 0<\norm{f}^2\tilde{\lambda}_0\leq \norm{D_rf}^2+\norm{D_r\widehat{f}}^2\leq \norm{f}^2\tilde{\lambda}_{N-1},\label{ineq1}\end{align}
where $0<\tilde{\lambda}_0\leq \tilde{\lambda}_1\leq...\leq \tilde{\lambda}_{N-1}$ are the ordered real eigenvalues of $L+\Delta$.
Furthermore, the bounds are sharp.
\end{thm}
\PROOF Noting that
\begin{align*}\norm{D_rf}^2&=\inne{D_rf,D_rf}=\inne{f,\chi \Delta \chi^*f}=\inne{\widehat{f},\Delta \widehat{f}}
\mbox{ and that }\norm{D_r\widehat{f}}^2=\inne{\widehat{f},L\widehat{f}},\end{align*}
 we have
\begin{align*}\norm{D_rf}^2+\norm{D_r\widehat{f}}^2=
\inne{\widehat{f},(L+\Delta)\widehat{f}}.\end{align*}
Assuming $\tilde{\lambda}_0>0$, Inequality (\ref{ineq1}), and 
its sharpness, follow directly 
from $L+\Delta$ being symmetric and positive semidefinite, and by
 applying the properties of the Rayleigh quotient to $L+\Delta$.
To prove positivity of $\tilde{\lambda}_0$, note that for $\inne{\widehat{f},(L+\Delta)\widehat{f}}=0$ we must
have $\inne{h,\Delta h}=0=\inne{h,L h}$ for some $h\neq 0$. This is impossible as we have, for
non-zero $h$, $\inne{h,\Delta h}=0$ if and only if  $h=c[1,0,....,0]^*$ for some $c\neq 0$. 
This implies $\inne{h,L h}=deg(v_0)c^2> 0$ due to the
connectivity of the graph.\hfill $\blacksquare$


Alternatively, if we consider the normalized Laplacian $\mathcal{L}$, we define 
a slightly different notion of the difference operator in order
to reflect the normalized structure when using the normalized Laplacian.
For a function $f\in l^2(G)$, define the normalized graph
difference operator as $$D_{nr}=D^{-1/2}D_r=D^{-1/2}W^{1/2}M.$$ 
We can then prove (in Theorem \ref{Koprowski2}) a graph
differential normalized Fourier transform inequality of the form of Theorem \ref{Koprowski1}.

\begin{thm}\label{Koprowski2}
Let $G$ be a simple, connected, and undirected graph. Then, for any non-zero function $f\in l^2 (G)$, the following inequalities hold:
\begin{align} 0<\norm{f}^2\tilde{\mu}_0\leq \norm{D_{nr}f}^2+\norm{D_{nr}\overset{*}{f}}^2\leq \norm{f}^2\tilde{\mu}_{N-1},\label{ineq2}\end{align}
where $0<\tilde{\mu}_0\leq \tilde{\mu}_1\leq...\leq \tilde{\mu}_{N-1}$ are the ordered real eigenvalues of $\mathcal{L}+\mathcal{D}$.
Furthermore, the bounds are sharp.
\end{thm}
\PROOF After suitable substitution of normalized elements, the result follows 
from arguments similar to the proof of Theorem \ref{Koprowski1}.\hfill $\blacksquare$

\section{Graph Frame Uncertainty Principles}\label{GFUP}

As a generalization of the work by Lammers and Maeser
in \cite{lammers2011uncertainty}, we show that the modified Laplacian operator  $L+\Delta$ will dictate an additive uncertainty
principle for frames.
Let \begin{align*}E=\EE\end{align*}
be a $d\times N$ matrix whose columns form a Parseval frame for $\complexes^d,$ i.e., $EE^*=I_{d\times d}$.
Let $T$ be the $N\times N$ permutation matrix, viz., 
$$T=\cT.$$
Let $\mathcal{D}=T^0-T$, then $\mathcal{D}^*=T^0-T^{N-1},$ and 
the classical Laplacian in the discrete setting is given by $L_c=\mathcal{D}^*\mathcal{D}=2T^0-T-T^{N-1}$. Let $\Delta_c$ be the 
diagonalization of $L_c$.
Let $\norm{\cdot}_{fr}$ denote the Frobenius norm. Let $DFT$ denote the unitary discrete Fourier transform
matrix. The following result holds.
\begin{thm}(Lammers and Maeser  \cite{lammers2011uncertainty})
For fixed dimension $d$ and $N\geq d\geq 2$, the following inequalities hold for all $d\times N$ Parseval frames:
\begin{align}0<G(N,d)&\leq \norm{\mathcal{D}DFT E^*}_{fr}^2+\norm{\mathcal{D} E^*}_{fr}^2\notag\\&\leq H(N,d)\\& \leq 8d\notag.\end{align}
Furthermore, the minimum (maximum)
occurs when columns of $E^*$ are the $d$ orthonormal eigenvectors corresponding to the $d$ smallest (largest) eigenvalues of $L_c+\Delta_c$. 
The constant $G(N,d)$ is the sum of those $d$ smallest eigenvalues, and $H(N,d)
$ is the sum of those $d$ largest eigenvalues. \label{thm2}
\end{thm}




To extend the inequalities in Theorem \ref{thm2} to the graph Fourier transform setting, 
we apply $D_r$ to the frame's conjugate transpose $E^*$ and to the graph 
Fourier transform $\chi^* E^*$, and then find bounds for the Frobenius norms.
\begin{thm}\label{graph_frame_up}
For any graph $G$ as in Theorem \ref{Koprowski1}, the following inequalities hold for all  $d\times N$ Parseval frames $E$:
\begin{equation}
\sum_{j=0}^{d-1}\tilde{\lambda}_j\leq \norm{D_r\chi^* E^*}^2_{fr}+\norm{D_rE^*}^2_{fr}\leq \sum_{j=N-d}^{N-1}\tilde{\lambda}_j\label{gup}
,\end{equation} where $\set{\tilde{\lambda_j}}$ is the ordered set of real, positive eigenvalues of $L+\Delta$. Furthermore, these
bounds are sharp.
\end{thm}
\PROOF Using the trace formulation of the Frobenius norms yields
\begin{align}
\norm{D_r\chi^* E^*}^2_{fr}+\norm{D_rE^*}^2_{fr}&=tr(E\chi D_r^*D_r \chi^* E^*)
+tr(D_rE^*ED_r^*).\end{align}
Using the invariance of the trace when reordering products, we have
\begin{align*}\norm{D_r\chi^* E^*}^2_{fr}+\norm{D_rE^*}^2_{fr}
&=tr(L\chi^* E^*E\chi)+tr(LE^*E)\\
&=tr(L \chi^*E^*E\chi)+ tr(\chi\Delta\chi^* E^*E)\\
&=tr((L+\Delta)\chi^*E^*E\chi).\end{align*}
The operator $\Delta+L$ is real, symmetric, and positive semidefinite.
By the spectral theorem, it has an orthonormal eigenbasis $P$ that, upon conjugation, diagonalizes $\Delta+L$: $$P^*(\Delta+L)P=\tilde{\Delta}=
\mbox{diag}( \tilde{\lambda_0},\tilde{\lambda}_1,..., \tilde{\lambda}_{N-1}).$$
Hence, we have
\begin{align*}\norm{D_r\chi E^*}^2_{fr}+\norm{D_rE^*}^2_{fr}&=tr((\Delta+L)\chi^*E^*E\chi)
=tr(P\tilde{\Delta} P^*\chi^*E^*E\chi)\\
&=tr(\tilde{\Delta} P^*\chi^*E^*E\chi P)
=\sum_{j=0}^{N-1} \pa{K^*K}_{j,j}\tilde{\lambda}_j,
\end{align*} where $K=E\chi P$.
The matrix $K$ is a Parseval frame because unitary transformations of Parseval
frames are Parseval frames. Therefore, $tr(K^*K)=tr(KK^*)=d.$ $K^*K$ is also the product
of matrices with operator norm $\leq 1$. Therefore, each of the entries, $\pa{K^*K}_{j,j},$ satisfies
$0\leq \pa{K^*K}_{j,j}\leq 1$.
Hence, minimizing (maximizing) $\sum_{j=0}^{N-1} \pa{K^*K}_{j,j}\tilde{\lambda}_j$
is achieved if $$\pa{K^*K}_{j,j}=\begin{cases} 1 & j<d ~(j\geq N-d)\\0&j\geq d~(j<N-d).\end{cases}$$
Choosing $E$ to be the first (last) $d$ rows of $(\chi P)^*$ accomplishes this. The positivity of the bounds
follows from the proof of Theorem \ref{Koprowski1}\hfill $\blacksquare$

A similar result holds for the normalized graph Laplacian.

\begin{thm}\label{graph_frame_up_norm}
For any graph $G$ as in Theorem \ref{Koprowski1}, the following inequalities hold for all  $d\times N$ Parseval frames $E$:
\begin{equation}
\sum_{j=0}^{d-1}\tilde{\mu}_j\leq \norm{D_{nr}\mathcal{F}^* E^*}^2_{fr}+\norm{D_{nr}E^*}^2_{fr}\leq \sum_{j=N-d}^{N-1}\tilde{\mu}_j\label{gup2}
,\end{equation} where $\set{\tilde{\mu_j}}$ is the ordered set of real, positive eigenvalues of $\mathcal{L}+\mathcal{D}$. Furthermore, these
bounds are sharp.
\end{thm}

\PROOF After suitable substitution of normalized elements, the result follows 
from arguments similar to the proof of Theorem \ref{graph_frame_up}.\hfill $\blacksquare$

\section{Feasibility Region}\label{feasibility}

We extend the concept of the feasibility region for graph and spectral spreads
from Agaskar and Lu \cite{agaskar2013spectral}. As such, the structure of the arguments used
here follows the structure of Agaskar and Lu's proofs.
Define the \textit{difference operator feasibility region} $FR$ as follows: $$FR=\set{(x,y): 
\norm{D_rf}^2=x\mbox{ and }\norm{D_r\widehat{f}}^2=y\mbox{ for some unit normed }f\neq0\in l^2(G)}.$$
Our analysis relies on
a key lemma, which in turn, relies on the following theorem due to Barvinok \cite{AB}.
\begin{thm}(Barvinok \cite{AB})\label{Barv}
Let $\mbox{Sym}_N$ be the set of real $N\times N$ symmetric matrices and
let $S_{+}^N$ be the subset of positive semidefinite symmetric matrices.
Suppose that $R> 0$ and $N\geq R+2$. Let $\mathcal{H}\subset \mbox{Sym}_N $
be an affine subspace such that codim$(\mathcal{H})= \binom{R+2}{2}$. 
If $S_{+}^N\cap \mathcal{H}$ is nonempty and bounded, 
then there exists a matrix $M\in S_{+}^N\cap \mathcal{H}$ of rank less than or equal to $R$.
\end{thm}

\begin{lem}\label{Keylem}
 Let $M_1$ and $M_2$ be rank one positive semidefinite matrices such that 
 \begin{align}
  tr(M_i)=1,~ tr(\Delta M_i)=x_i \mbox{ and } tr(LM_i)=y_i \mbox{ for } i=1,2.\label{keylem}
 \end{align}
 Then, for any $\beta\in[0,1]$, there exists a rank one positive semidefinite matrix $M$
 satisfying
 \begin{align}
  tr(M)=1,~ tr(\Delta M)=x \mbox{ and } tr(LM)=y, \label{keylem2}
 \end{align}
where $x=\beta x_1+(1-\beta)x_2$ and $y=\beta y_1 + (1-\beta)y_2.$
\end{lem}

\PROOF Let each positive semidefinite matrix $M_i$ satisfy
(\ref{keylem}).
For any $\beta\in [0,1]$, let $M'=\beta M_1+(1-\beta)M_2.$ 
Clearly, $M'\in S_+^N$ by the convexity of $S_+^N$ and if we 
let $x=\beta x_1+(1-\beta)x_2$ and $y=\beta y_1+(1-\beta)y_2$ 
then $$M'\in \mathcal{H}=\set{M\in\mbox{Sym}_N: tr(M)=1,tr(\Delta M)=x,\mbox{ and } tr(LM)=y}.$$
 By the linear independence of $I$, $\mathcal{L}$, and $\Delta$, we have that $\mathcal{H}$ is an affine subspace 
of $\mbox{Sym}_N$ with codimension 3. 
Hence, we have that $S_{+}^N\cap \mathcal{H}\neq \emptyset.$ 
Noting that any element of $S_{+}^N\cap \mathcal{H}$
has nonnegative eigenvalues which must sum to 1,
the boundedness of this subspace is straightforward to show:
$$\forall M\in S_{+}^N\cap \mathcal{H},~ \norm{M}^2_{fr}=tr(M^2)\leq tr(M)=1.$$
By Theorem \ref{Barv}, we conclude that there exists a matrix $M$ of rank one
satisfying (\ref{keylem2}). \hfill $\blacksquare$

 The importance of Lemma \ref{Keylem} to our arguments shall lie in the fact that
 any matrix $M$ satisfying (\ref{keylem2}) has eigenvalue $1$ with multiplicity
 $1$ and eigenvalue $0$ with multiplicity $N-1$.
 Hence, if $Mg=g$ for a unit vector $g$, then $M=gg^*$.  
We shall prove some key properties of the difference operator feasibility region.
\begin{thm} \label{DNFR} Let $FR$ be the difference operator feasibility region for a simple, connected graph
$G$ with $N$ vertices.
Then, the following properties hold.
\begin{itemize}\item[a)] $FR$ is a closed subset of $[0, \lambda_{N-1}]\times[0, \lambda_{N-1}]$ where
$\lambda_{N-1}$ is the maximal eigenvalue of the Laplacian $L$.

\item[b)] $y=0$ and $x=\frac{1}{N}\sum_{j=0}^{N-1}\lambda_j$ is the only point on the horizontal axis in $FR$. $x=0$ and $y=L_{0,0}$
 is the only point on the vertical axis in $FR$.

\item[c)] $FR$ is in the half plane defined by $x+y\geq \tilde{\lambda}_0>0$
with equality if and only if $\widehat{f}$ is in the eigenspace 
associated with $\tilde{\lambda}_0$.
\item[d)] If $N=2$, then $FR$ is the circle
given by $$FR=\set{(\alpha(a-b)^2, 2\alpha b^2): a^2+b^2=1\mbox{ and } \alpha>0}.$$
If $N\geq 3,$ then $FR$ is a convex region. 
\end{itemize}
\end{thm}

\PROOF Recall that $\norm{D_r f}^2=\inne{f,Lf}=\inne{\widehat{f},\Delta \widehat{f}},$ 
and that $\norm{D_r\widehat{f}}^2=\inne{\widehat{f}, L\widehat{f}}.$
Note that the operation $f\mapsto \widehat f$ is an isomorphism
of the unit sphere in $l^2(G)$. Hence, for the entirety of this proof
we rely on the fact that if a unit normed $g\in l^2(G)$ (respectively, 
a unique unit normed $g\in l^2(G)$) achieves a value
in the feasibility region for $\inne{g,\Delta g}$, and
for $\inne{g,L g},$ then there exists a
unit normed  $f\in l^2(G)$ (respectively, a unique 
unit normed  $f\in l^2(G)$)
that achieves the same values for $\norm{D_rf}^2$ and $\norm{D_r\widehat{f}}^2$ respectively,
viz., $f=\chi g$.

\begin{itemize}\item[a)]By the properties of the Rayleigh quotient, we have
for any unit normed $g\in l^2(G)$ that
$$0=\lambda_0\leq \inne{g,\Delta g}\leq \lambda_{N-1}.$$
The maximum is attained if $g=[0,...,0,1]'$. Similarly,
we have that
$$0=\lambda_0\leq \inne{g,Lg}\leq \lambda_{N-1},$$
and the maximum is attained if $g$ is in the eigenspace associated
with $\lambda_{N-1}$ for $L$.
Hence, $FR\subset [0, \lambda_{N-1}]\times[0, \lambda_{N-1}]$.
$FR$ is closed because it is the image of a continuous mapping from
the closed unit sphere of $l^2(G)$ into $\reals^2$.

\item[b)] $\inne{g,Lg}=0$ if and only if $g=\pm\frac{1}{\sqrt{N}}[1,...,1]'.$
Hence, we have $$x=\inne{g,\Delta g}=\frac{1}{N}[1,...,1]
\begin{bmatrix} \lambda_0 \\
\lambda_1\\
\vdots\\
\lambda_{N-1}\end{bmatrix}=\frac{1}{N}\sum_{j=0}^{N-1}\lambda_j.$$
$x=0$ if and only if $g=\pm[1,0,...,0]'$. Hence, we have 
$y=L_{0,0},$ which is the degree of the first vertex of $G$. 

\item[c)] This follows directly from Theorem \ref{Koprowski1}.

\item[d)] If $N=2$, then the only simple connected graph is 
a graph with one edge of weight $\alpha>0$. Applying the corresponding operators to
any unit normed vector $g=[a,b]^*$ yields the desired result.

To show that $FR$ is convex for $N\geq 3$, 
we shall formulate this as a problem in $\mbox{Sym}_N.$
First, we note that showing convexity is equivalent to showing the following:
if $g_1,g_2\in l^2(G)$ and 
\begin{align}
\inne{g_i,g_i}=1,~\inne{g_i,\Delta g_i}=x_i,\mbox{ and } \inne{g_i, Lg_i}=y_i\mbox{ for } i=1,2,\label{convec4}
\end{align}
then for any $\beta\in[0,1]$, we can always find a function $g\in l^2(G)$ satisfying
\begin{align}
\inne{g,g}=1,~\inne{g, \Delta g}=x,\mbox{ and } \inne{g, Lg}=y\label{convec3}
\end{align}
where $x=\beta x_1+(1-\beta)x_2$ and $y=\beta y_1 + (1-\beta)y_2.$ 
Let $g_1$ and $g_2$ satisfy (\ref{convec4}) and set $M_1=g_1g_1^*$ and $M_2=g_2g_2^*.$
Applying Lemma \ref{Keylem} to $M_1$ and $M_2$ yields $M=gg^*$ where
$M$ satisfies (\ref{keylem2}). 
This is easily shown, via the cyclic properties of the trace operator,
to be equivalent to $g$ satisfying (\ref{convec3}).
$\phantom{alkhd}$\hfill $\blacksquare$
\end{itemize}

\begin{figure}
 \begin{tikzpicture}[scale=1.5]
    \coordinate (y) at (0,5);
    \coordinate (x) at (6.2,0);
    \draw[<->] (y) node[above] {$\langle g, L g\rangle$} -- (0,0) --  (x) node[right]
    {$\langle g, \Delta g\rangle$};

    \path
    coordinate (start) at (0,4)
    coordinate (c1) at +(1,.2)
    coordinate (c2) at +(2.5,.1)
    coordinate (slut) at (4,0)
    coordinate (c3) at (4.5, .05)
    coordinate (c4) at (5.5, .1)
    coordinate (top) at (6,1);

    \draw[important line, red, thick] (start) .. controls (c1) and (c2) .. (slut);
    \draw[important line, red, thick] (slut) .. controls (c3) and (c4) .. (top);
    \draw[thick, dashed] (1.5,0) -- (1.5,5);
    \draw[thick, dashed] (6,0) -- (6,5);
    \draw[->] (5, 3) node[above] {$\omega(x)$} -- (3,.1);
    \filldraw (1.5, 0) circle (2pt) node[below] {$x_0$};
    \filldraw (6, 0) circle (2pt) node[below] {$\lambda_{N-1}$};
    \filldraw (4,0) circle (2pt) node[align=center, below] {$(\sum_j\frac{\lambda_j}{N},0)$};
    \filldraw (0,4) circle (2pt) node[left] {$(0, L_{0,0})$};
    \filldraw (1.5, .9) circle (2pt) node[right] {$(x_0, \langle g', L g'\rangle)$};
    \filldraw (1.5, 3) circle (2pt) node[right] {$(x_0, \langle g, L g\rangle)$};
   
    \end{tikzpicture}\caption{The DUC (red) for a simple and connected graph $G$}
\label{DUCpic}
\end{figure}
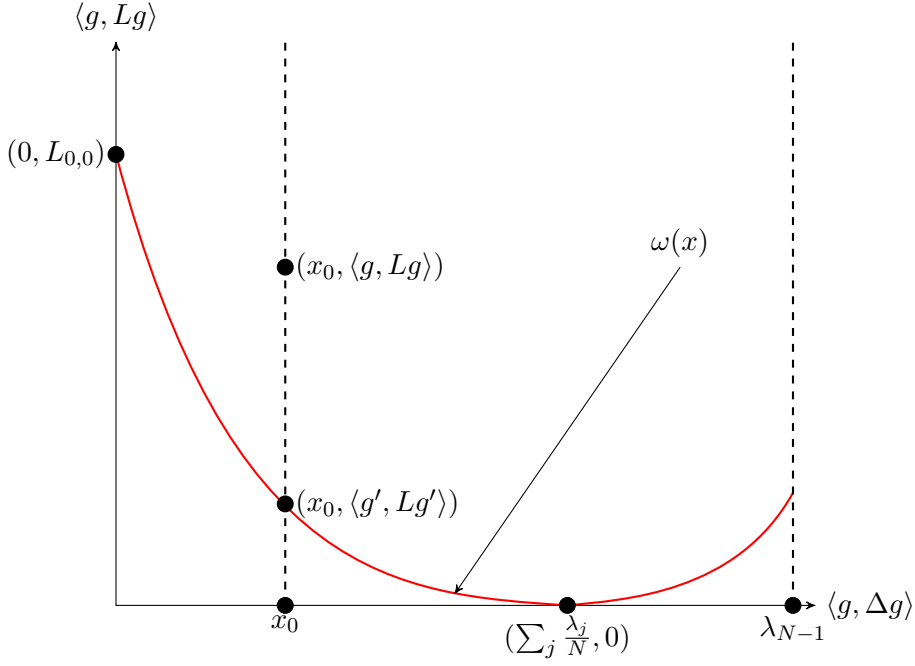

We now turn our attention to the lower boundary of $FR$.
The differential uncertainty curve (DUC)
$\omega(x)$ is defined as follows:
\begin{align*} \forall x\in [0,\lambda_{N-1}],~ \omega(x)=\inf_{g\in l^2(G)}\inne{g,Lg}
\mbox{ subject to } \inne{g,\Delta g}=x\mbox{ and }\norm{g}=1.\end{align*} See Figure (\ref{DUCpic})
for a sample uncertainty curve.
Given a fixed $x\in [0,\lambda_{N-1}]$, we say  $g'$ \textit{attains} the DUC
if for all $g$ with $\inne{g,\Delta g}=x$ we have
$\inne{g',Lg'}\leq\inne{g,Lg}.$
We shall show that for all $x\in[0,\lambda_{N-1}]$, there exists a 
function attaining the DUC.
In fact, we shall show 
that certain eigenfunctions of the matrix valued
function $K(\alpha)=L-\alpha \Delta$
will attain the DUC for every
value of $x$. Hence, we shall show that for all $x\in (0,\lambda_{N-1}),$
$$\omega(x)=\min_{g\in l^2(G)} \inne{g,Lg}\mbox{ subject to }
\inne{g, \Delta g}=x\mbox{ and } \inne{g,g}=1.$$

We begin classifying $\omega(x)$ by motivating the use of the
operator $K(\alpha)$. Finding values that attain the differential
uncertainty curve amounts to solving a quadratically constrained
convex optimization problem. We achieve this by defining the
following Lagrangian function, and setting its gradient equal to
zero. Define the DUC \textit{Lagrangian}, $\Gamma$, as  
$$\Gamma(g,\alpha,\beta)=\inne{g,Lg}-\alpha(\inne{g,\Delta g}-x)-
\beta(\inne{g,g}-1).$$ Upon taking the gradient with respect to $g$
and setting the gradient equal to zero, we have for some optimal $g'$ 
that
$$\nabla_g(\Gamma(g,\alpha, \beta))(g')
=2Lg'-2\alpha\Delta g'-2\beta g'=0$$
and $$K(\alpha)g'=(L-\alpha \Delta)g'=\beta g'.$$
Thus, the minimizer of the quadratically constrained problem
is an eigenfunction of the operator $K(\alpha)$.
Define $m(\alpha)$ to be the minimal eigenvalue of $K(\alpha)$,
and define $\sigma(\alpha)$ to be its associated eigenspace.
We shall prove that a function $g$ attains the DUC if and only
if it in $\sigma(\alpha)$. In order to prove this, we
shall rely heavily on analysis of the functions
$H_+$ and $H_-$ defined as follows:
\begin{align}H_+(\alpha)=\max_{g\in \sigma(\alpha):\norm{g}=1}\inne{g,\Delta g}
\mbox{ and } H_-(\alpha)=\min_{g\in \sigma(\alpha):\norm{g}=1}\inne{g,\Delta g},\label{HHH}\end{align}
which measure the maximal, respectively, minimal values
that can be achieved by eigenfunctions in $\sigma(\alpha)$.

\begin{lem}\label{HHlem}
The following properties hold for $H_+(\alpha)$ and $H_-(\alpha)$.
\begin{itemize}\item[a)] For all $\alpha\in \reals,$ $H_+(\alpha)$
and $H_-(\alpha)$ are increasing functions.
\item[b)] $\lim_{\alpha\to\infty} H_{\pm}(\alpha)=\lambda_{N-1}$, and 
$\lim_{\alpha\to -\infty} H_{\pm}(\alpha)=0$.
\item[c)] On any finite interval $[a,b]$, 
the functions $H_\pm$ differ on at most a finite number of points 
denoted by $\Sigma=\set{b_1,...,b_k}$ for some $k\geq 0$.
For all $\alpha\not\in\Sigma$, the following holds: $H_+(\alpha)=H_-(\alpha)=-m'(\alpha)$.
\end{itemize}
\end{lem}
\PROOF \begin{itemize}
\item[a)] For $\alpha_1<\alpha_2$, we take any $\nu_1\in \sigma(\alpha_1)$ 
and $\nu_2\in \sigma(\alpha_2)$,
and we have, by the Rayleigh quotient for symmetric matrices, that
$$\inne{\nu_2,K(\alpha_1)\nu_2}\geq m(\alpha_1)=\inne{\nu_1,K(\alpha_1)\nu_1}.$$
Similarly, we have
$$-\inne{\nu_2,K(\alpha_2)\nu_2}=- m(\alpha_2)\geq-\inne{\nu_1,K(\alpha_2)\nu_1}.$$
Combining the inequalities yields
\begin{align}\inne{\nu_2,(K(\alpha_1)-K(\alpha_2))\nu_2}\geq \inne{\nu_1,(K(\alpha_1)-K(\alpha_2))\nu_1}.\label{52}\end{align}
Noting that $K(\alpha_1)-K(\alpha_2)=\pa{\alpha_2-\alpha_1}\Delta$, and substituting into (\ref{52}) yields
\begin{align}\inne{\nu_2,\Delta\nu_2}\geq \inne{\nu_1,\Delta\nu_1}\notag\end{align}
Upon specializing to the unit norm eigenfunctions that attain the maximization in (\ref{HHH})  we have
\begin{align} H_+(\alpha_2)=\inne{\nu_2,\Delta\nu_2}\geq \inne{\nu_1,\Delta\nu_1}=H_+(\alpha_1)\notag\end{align}
Similarly, upon specializing to the unit norm eigenfunctions that attain the minimum in (\ref{HHH})
we have that $H_-(\alpha_2)=\inne{\nu_2, \Delta \nu_2}\geq\inne{\nu_1,\Delta\nu_1}= H_-(\alpha_1).$

\item[b)] Let $\alpha\in\reals$, then we clearly have 
$$H_+(\alpha)\geq H_-(\alpha)\geq 0$$ by the positive 
semidefinite property of $\Delta$.
Let $\nu\in \sigma(\alpha)$ be unit normed. 
Recall that the canonical first basis vector $e_0$
spans the null space of $\Delta$ 
and hence $\inne{e_0,\Delta e_0}=0$. 
For any unit norm $\nu\in\sigma(\alpha),$
we have $\inne{\nu,L\nu}\geq 0$,
and if $\alpha<0$, we have $-\alpha\inne{\nu,\Delta \nu}\geq 0$.
Thus by the properties the Rayleigh quotient we have
\begin{align} 0\leq -\alpha\inne{\nu,\Delta \nu}
\leq \inne{\nu, K(\alpha)\nu}\leq 
\inne{e_0, K(\alpha)e_0}
=L_{o,o}+0=L_{o,o}.\label{l00}\end{align} 

Multiplying (\ref{l00}) by $-\frac{1}{\alpha}$ yields
\begin{align*} 0\leq \inne{\nu,\Delta\nu}&\leq 
-\frac{1}{\alpha}L_{0,0}.\end{align*}
Since this is valid for all $\nu\in\sigma(\alpha)$
we have 
\begin{align*}0\leq H_-(\alpha)\leq H_+(\alpha)\leq 
-\frac{L_{0,0}}{\alpha}\end{align*}
As $\alpha\to-\infty$, we squeeze $H_\pm(\alpha)$ to zero as desired.

For the limit as $\alpha\to\infty$, recall that the last
canonical eigenfunction $e_{N-1}$ is in the
eigenspace of $\lambda_{N-1}$ for $\Delta.$ 
Hence, $\inne{e_{N-1},\Delta e_{N-1}}=\lambda_{N-1}$, 
and we have
\begin{align*}\inne{\nu, K(\alpha)\nu}&\leq \inne{e_{N-1}, 
K(\alpha)e_{N-1}} \\
&= \inne{e_{N-1},  Le_{N-1}}-\alpha\lambda_{N-1}\\
&= L_{N-1,N-1}-\alpha\lambda_{N-1}.
\end{align*}
Adding $\pa{\alpha\inne{\nu, \Delta\nu}-L_{N-1,N-1}}$
to both sides yields
\begin{align}\inne{\nu, L\nu}-L_{N-1,N-1}\leq \alpha
\pa{\inne{\nu,\Delta\nu}-\lambda_{N-1}}
\leq 0\label{negative}\end{align} 
where the last inequality in (\ref{negative}) is due to $\alpha>0$ and 
the properties of the Rayleigh quotient. Taking
the absolute value of both sides, and dividing by $\alpha$
yields $$\abs{\frac{\inne{\nu, L\nu}-L_{N,N}}{\alpha}}
\geq \abs{\inne{\nu,\Delta\nu}-\lambda_{N-1}}\geq 0.$$
The desired result follows by taking $\alpha\to\infty$.

\item[c)] We use eigenvalue perturbation results  such as those in \cite{lancaster} to establish that $m(\alpha)$ is analytic for
$[a,b]\cap(\Upsilon)^c$ where $\Upsilon$ is a finite subset of $[a,b]$. 
$K(\alpha)$ is real, is linear in $\alpha$, 
hence analytic, and it is symmetric. 
By Theorem 2 on page 404 of \cite{lancaster}, 
there exist $N$ analytic functions
$\xi_0(\cdot),...,\xi_{N-1}(\cdot)$ 
and $N$ analytic vector valued functions
$w_0(\cdot),..., w_{N-1}(\cdot)$ such that 
\begin{align} \forall \alpha\in\reals~K(\alpha)w_j(\alpha)=\xi_j(\alpha)w_j(\alpha)\end{align}
and $$\inne{w_j(\alpha),w_k(\alpha)}=
\begin{cases} 0 &\mbox{ if } j\neq k\\ 1 &\mbox{ if } j=k\end{cases}.$$

Let $[a,b]$ be an arbitrary finite interval in $\reals$, and fix
$\alpha_0\in(a,b)$. If $\sigma(\alpha_0)$ is one dimensional,
then exactly one eigenvalue function $\xi_j(\alpha_0)$ equals 
$m(\alpha_0)$. By the analycity of all the eigenvalue
functions, there exists some $\delta$ ball about $\alpha_0$, such
that if $\abs{\alpha-\alpha_0}<\delta$ then 
$\xi_j(\alpha)<\xi_k(\alpha)$ for $k\neq j$.
Hence, $m(\alpha)=\xi_j(\alpha)$ for $\alpha\in(\alpha_0-\delta,
\alpha_0+\delta)$ and therefore $m(\alpha)$ is analytic on the $\delta$ ball.

If $\sigma(\alpha_0)$ has dimension greater than one, then more than
one eigenvalue function from $\xi_l(\alpha_0)$ for $l=0,...,N-1$ attains the
value $m(\alpha_0)$. In this case, $m(\alpha)$ may not be analytic in
any neighborhood of $\alpha_0$. For instance, if two of the eigenvalue
functions cross at exactly $\alpha_0$, then there is no derivative
for $m(\alpha)$ at $\alpha_0$. Define $\rho_j(\alpha)$ for 
$j=0,...,d\leq N-1$ as the distinct eigenvalue functions of
$K(\alpha)$, and let $n_j$ be the multiplicity of each function.
For $[a,b]\subset\reals$, define
$$\Upsilon=\bigcup_{0\leq i<j\leq d}\set{\alpha\in[a,b]:\rho_i(\alpha)
=\rho_j(\alpha)}.$$
As defined, $\Upsilon$ has finite order. Indeed, if $\abs{\Upsilon}=
\infty$ then at least two of the $\rho_j$'s would be equal on an
infinite set of points on the interval, and therefore would be equal
on the interval because both are analytic.

To conclude the proof, we shall relate $m(\alpha)$ to
$H_\pm(\alpha)$. For fixed $\alpha_0\in [a,b]$, we
without loss of generality, assume the first $k+1$ distinct
eigenvalue functions $\rho_i$ for $i=0,...,k$ intersect at $\alpha_0$,
and are minimal. That is to say, $\rho_i(\alpha_0)=m(\alpha_0)$.
The associate eigenfunction functions are denoted by 
$w_{i,j}(\alpha)$ for $i=0,...,k$, and $j=1,...,n_i$.
Hence, $w_{ij}(\alpha_0)$ form an orthonormal basis for
$\sigma(\alpha_0)$. 
Therefore, if $\nu\in\sigma(\alpha_0)$ is unit normed,
we have $$\nu=\sum_{i=0}^k\sum_{j=1}^{n_i}c_{ij}w_{ij}(\alpha_0).$$
The coefficients are $c_{ij}=\inne{\nu,w_{ij}(\alpha)},$
and, therefore, we have $\sum_{i=0}^k\sum_{j=1}^{n_i}c_{ij}^2=\norm{\nu}^2=1.$
We define the analytic function $\nu(\alpha)$ such that 
$\nu(\alpha_0)=\nu$ as follows:
$$\nu(\alpha)=\sum_{i=0}^k\sum_{j=1}^{n_i}c_{ij}w_{ij}(\alpha).$$
Applying $K(\alpha)$ to $\nu(\alpha)$ yields
\begin{align}K(\alpha)\nu(\alpha)&=
\sum_{i=0}^k\sum_{j=1}^{n_i}c_{ij}K(\alpha)w_{ij}(\alpha)
=\sum_{i=0}^k\sum_{j=1}^{n_i}c_{ij}\rho_i(\alpha)w_{ij}(\alpha).
\label{Knu}\end{align}
We apply the product rule to differentiate equation (\ref{Knu})
which yields
\begin{align} K'(\alpha)\nu(\alpha)+K(\alpha)\nu'(\alpha)=
\sum_{i=0}^k\sum_{j=1}^{n_i}c_{ij}\rho_{i}'(\alpha)w_{ij}(\alpha)+
\sum_{i=0}^k\sum_{j=1}^{n_i}c_{ij}\rho_{i}(\alpha)w'_{ij}(\alpha).
\label{Knuprime}\end{align}
Evaluating equation (\ref{Knuprime}) at $\alpha_0$ yields
\begin{align} K'(\alpha_0)\nu(\alpha_0)+K(\alpha_0)\nu'(\alpha_0)=
\sum_{i=0}^k\sum_{j=1}^{n_i}c_{ij}\rho_{i}'(\alpha_0)w_{ij}(\alpha_0)+
\sum_{i=0}^k\sum_{j=1}^{n_i}c_{ij}\rho_{i}(\alpha_0)w'_{ij}(\alpha_0).
\label{Knuprime0}\end{align}
Noting that $K'(\alpha)=-\Delta$, $\rho_i(\alpha_0)=m(\alpha_0)$, 
$$\inne{\nu(\alpha_0),K(\alpha_0)\nu'(\alpha_0)}=
\inne{K(\alpha_0)\nu(\alpha_0),\nu'(\alpha_0)}=
m(\alpha_0)\inne{\nu(\alpha_0),\nu'(\alpha_0)},$$
and that
$$\inne{\nu(\alpha_0),
\sum_{i=0}^k\sum_{j=1}^{n_i}c_{ij}\rho_{i}'(\alpha_0)w_{ij}(\alpha_0)}=
\sum_{i=0}^k\sum_{j=1}^{n_i}c^2_{ij}\rho_{i}'(\alpha_0),$$
we have that the inner
product of $\nu(\alpha_0)$ with the left and right hand sides of 
equation (\ref{Knuprime0}) yields

\begin{align}&-\inne{\nu(\alpha_0),\Delta\nu(\alpha_0)}+m(\alpha_0)
\inne{\nu(\alpha_0),\nu'(\alpha_0)}\notag\\
&=\sum_{i=0}^k\sum_{j=1}^{n_i}c^2_{ij}\rho_{i}'(\alpha_0)
+m(\alpha_0)\sum_{i=0}^{k}\sum_{j=1}^{n_i}c_{ij}\inne{\nu(\alpha_0),
w_{ij}'(\alpha_0)}.\label{summands}\end{align}

The second summands on the LHS and RHS of equation (\ref{summands})
are equal, hence we have 
$$\inne{\nu(\alpha_0),\Delta\nu(\alpha_0)}=
-\sum_{i=0}^k\sum_{j=1}^{n_i}c^2_{ij}\rho_{i}'(\alpha_0).$$
Since $\nu(\alpha_0)=\nu\in\sigma(\alpha_0)$ was arbitrary and 
unit normed, and since the dimension of $\sigma(\alpha_0)$ is 
finite, maximizing (respectively minimizing) $\inne{\nu,\Delta\nu}$
is achieved by maximizing (respectively minimizing) over
the $\rho_i(\alpha_0)$'s. Hence we have
$$H_+(\alpha_0)=\max_{0\leq i\leq k}-\rho_i'(\alpha_0),
\mbox{ and } H_-(\alpha_0)=\min_{0\leq i\leq k}-\rho_i'(\alpha_0).$$
Since all of the $\rho_i(\alpha)$ are distinct (except at
$\alpha_0$) in some neighborhood
$\mathcal{N}$ covering $\alpha_0$ and small enough
that $\mathcal{N}\cap\Upsilon=\alpha_0$ or $\emptyset$,
there exist $l,m\in\set{0,...,k}$ such that 
$$m(\alpha)=\begin{cases}\rho_l(\alpha) & \alpha\leq\alpha_0\\
\rho_m(\alpha) & \alpha\geq \alpha_0.\end{cases}$$
If for some $j\neq m$, $\rho_j'(\alpha_0)<\rho_m'(\alpha_0)$
then $\rho_j(\alpha)<\rho_m(\alpha)$ for some $\alpha\in\mathcal{N}
\cap[\alpha_0,\infty).$ This contradicts the fact that $m(\alpha)
=\rho_m(\alpha)$ on this interval. Similarly, if there exists some 
$j\neq l$, with $\rho_j'(\alpha_0)>\rho'_l(\alpha_0)$
there is a contradiction on $\mathcal{N}\cap (-\infty,\alpha_0].$
Hence, we have
$$H_+(\alpha)=-\rho'_m(\alpha)=-m'(\alpha)
\mbox{ for } \alpha\in\mathcal{N}\cap[\alpha_0,\infty),$$
and
$$H_-(\alpha)=-\rho'_l(\alpha)=-m'(\alpha)
\mbox{ for } \alpha\in\mathcal{N}\cap(-\infty,\alpha_0].$$
Right and left continuity follow from $\rho_m$ and $\rho_l$
having continuous derivatives. If $k=0$, that is, if only
one of the $\rho_i$ functions aligns with $m$ at $\alpha_0$,
then $m(\alpha)$ is analytic on $\mathcal{N}$ and we have
$H_-(\alpha)=H_+(\alpha)=-m'(\alpha)$ on $\mathcal{N}$.
If we denote $\Sigma$ as the set of $\alpha\in[a,b]$ for which
$H_-(\alpha)\neq H_+(\alpha)$, we must have $\Sigma\subseteq\Upsilon$
and therefore $\Sigma$ is a finite set.\hfill $\blacksquare$

\end{itemize}

We now prove that vectors in $\sigma(\alpha)$ characterize the DUC.

\begin{thm}\label{DUC} A unit normed function $f\in l^2(G)$ with $\norm{D_rf}^2=x
\in (0,\lambda_{N-1})$ achieves the DUC if and only if 
$\widehat{f}$ is a nonzero eigenfunction in $\sigma(\alpha)$ for some 
$\alpha\in\reals$.
\end{thm}
\PROOF As before, it suffices to show 
that a unit normed $\eta\in l^2(G)$ satisfying $\inne{\eta,\Delta 
\eta}=x\in (0,\lambda_{N-1})$ achieves the DUC if and only if 
$\eta\in\sigma(\alpha)$ for some $\alpha\in \reals$.

For the sufficient condition,
fix $\alpha\in\reals$. Then for any 
arbitrary unit norm $\eta\in l^2(G)$ we have
\begin{align} \inne{\eta,K(\alpha)\eta}&=\inne{\eta,L\eta}-\alpha 
\inne{\eta,\Delta\eta}.\notag\end{align} 
The Rayleigh quotient for $K(\alpha)$ is bounded sharply 
below by $m(\alpha)$. Hence we obtain
for any unit normed $\nu\in \sigma(\alpha)$,
\begin{align*}\inne{\nu, L\nu}-\alpha\inne{\nu,\Delta\nu}
=m(\alpha)\leq \inne{\eta,L\eta}-\alpha\inne{\eta,\Delta\eta}.\end{align*}
We assumed $\inne{\nu, \Delta\nu}=x$, and 
upon restricting $\eta$ to $\inne{\eta,\Delta\eta}=x,$ we have
$\inne{\eta,L\eta}\geq \inne{\nu, L\nu}.$
Hence, any unit normed $\nu\in\sigma(\alpha)$ achieves the DUC.

For the necessary condition, it suffices to show 
that for any function $\eta\in l^2(G)$ that 
achieves the DUC, there is an $\alpha$ 
and a unit norm $\nu\in \sigma(\alpha)$ such that 
$\inne{\nu,\Delta\nu}=\inne{\eta, \Delta\eta}=x$. Indeed,
having also assumed $\eta$ lies on the DUC, and being 
guaranteed that such a $\nu$ lies on the curve by the sufficient
condition, we have $\inne{\eta,L\eta}=\inne{\nu,L\nu}$, and 
hence
\begin{align*}\inne{\eta,K(\alpha)\eta}&=\inne{\eta,L\eta}-\alpha x=\inne{\nu, L\nu}-\alpha x
=\inne{\nu,K(\alpha)\nu}=q(\alpha).\end{align*}
Therefore, $\eta$ must also be a unit vector in $\sigma(\alpha)$.

We complete the proof by showing that for any $x\in(0,\lambda_{N-1})$ 
there is an $\alpha$ and a unit norm eigenfunction $\nu\in \sigma(\alpha)$
such that $\inne{\nu,\Delta\nu}=x$. 

Given $x\in(0,\lambda_{N-1})$, parts (b) and (c) of Lemma \ref{HHlem} 
ensure that there exist $a'<b'$ such that $H_-(a')<x<H_+(b')$ and that 
there are $a<b$ with $a'\leq a<b\leq b'$ such that
on the interval $[a,b]$ there exists at most one point $\beta\in[a,b]$ 
at which $H_-(\beta)<H_+(\beta)$. The interval $[H_-(a),H_+(b)]$ can 
be written as the union of three subintervals:
\begin{align*} [H_-(a),H_+(b)]=[H_-(a),H_-(\beta))\cup[H_-(\beta),H_+(\beta)]
\cup(H_+(\beta),H_+(b)].\end{align*}
Thus, $x$ must belong to one of these three intervals. If $x$ is in the 
first or third subinterval, the continuity of $H_-(\alpha)$ and $H_+(
\alpha)$, respectively, on these intervals guarantees for some 
$\alpha_-$, respectively, $\alpha_+$ that
$H_-(\alpha_-)=x$, respectively, $H_+(\alpha_+)=x$, on one of these intervals.
By the construction of the $H_{\pm}$ functions, this also guarantees a $\nu$
achieving the DUC exists. 

It remains to be shown that such an $\alpha$ and $\nu$ exist
for $x\in[H_-(\beta),H_+(\beta)]$. 
We set 
\begin{align*}\nu_+=\mbox{argmax}_{z\in \sigma(\beta),\norm{z}=1} \inne{z,\Delta z}
\mbox{ and }
\nu_-=\mbox{argmin}_{z\in \sigma(\beta),\norm{z}=1} \inne{z,\Delta z},\end{align*}
and define, for $\theta\in [0,\pi/2],$ the vector valued function,
\begin{align*}\nu(\theta)=\frac{\cos\theta \nu_++\sin\theta \nu_-}
{\pa{1+\sin(2\theta)\inne{\nu_+,\nu_-}}^{1/2}}.\end{align*}
The assumption that $H_-(\beta)\neq H_+(\beta)$ ensures that the 
denominator is nonzero. The numerator has norm squared given by
\begin{align*}\norm{\cos\theta \nu_++\sin\theta \nu_-}^2&=1+2\cos\theta\sin\theta\inne{\nu_+,\nu_-}
=1+\sin(2\theta)\inne{\nu_+,\nu_-},\end{align*} 
and so $\norm{\nu(\theta)}=1.$ 
Further, $\nu(\theta)$ is continuous and 
$\nu(\theta)\in \sigma(\beta)$. 
By continuity, the intermediate value theorem, and the fact 
that $\inne{\nu(\pi/2),\Delta \nu(\pi/2)}=H_-(\beta)$ 
and $\inne{\nu(0),\Delta \nu(0)}=H_+(\beta)$, 
we have that there exists $\theta_0\in[0,\pi/2]$ 
such that $\inne{\nu(\theta_0),\Delta \nu(\theta_0)}=x.$ 
\hfill$\blacksquare$

\section{The Complete Graph}\label{complete}
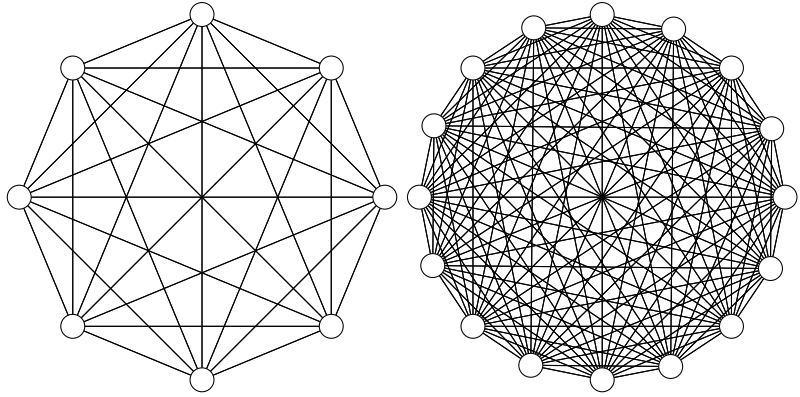
\begin{figure}[h]
\begin{center}\begin{tikzpicture}[transform shape, scale=.45]
  \foreach \number in {1,...,8}{
        \mycount=\number
        \advance\mycount by -1
  \multiply\mycount by 45
        \advance\mycount by 0
      \node[draw,circle,inner sep=0.25cm] (N-\number) at (\the\mycount:5.4cm) {};
    }
  \foreach \number in {1,...,7}{
        \mycount=\number
        \advance\mycount by 1
  \foreach \numbera in {\the\mycount,...,8}{
    \path (N-\number) edge[-] (N-\numbera)  edge[-] (N-\numbera);
  }
}
 
\end{tikzpicture} \begin{tikzpicture}[transform shape, scale=.45]
  \foreach \number in {1,...,8}{
        \mycount=\number
        \advance\mycount by -1
  \multiply\mycount by 45
        \advance\mycount by 0
      \node[draw,circle,inner sep=0.25cm] (N-\number) at (\the\mycount:5.4cm) {};
    }
  \foreach \number in {9,...,16}{
        \mycount=\number
        \advance\mycount by -1
  \multiply\mycount by 45
        \advance\mycount by 22.5
      \node[draw,circle,inner sep=0.25cm] (N-\number) at (\the\mycount:5.4cm) {};
    }
  \foreach \number in {1,...,8}{
        \mycount=\number
        \advance\mycount by 1
  \foreach \numbera in {\the\mycount,...,16}{
    \path (N-\number) edge[-] (N-\numbera)  edge[-] (N-\numbera);
  }
}
 \foreach \number in {9,...,15}{
        \mycount=\number
        \advance\mycount by 1
  \foreach \numbera in {\the\mycount,...,16}{
    \path (N-\number) edge[-] (N-\numbera)  edge[-] (N-\numbera);
  }
}
\end{tikzpicture}\end{center}
\caption{Unit weighted complete graphs with $N=8$ and $N=16$ vertices.}  \label{compgraph}
\end{figure}

Unit weighted graphs for which every vertex is connected directly to every other vertex,
as in Figure \ref{compgraph}, are referred to as \textit{complete graphs}.  A complete graph with $N$ vertices
has graph Laplacian $L=NI-O_{N\times N}$ where $O_{N\times N}$ is an $N\times N$ matrix each of whose elements is $1$. The minimal polynomial $m(x)$
for $L$ is given by $m(x)=x(x-N)$, and the characteristic polynomial is $c(x)=x(x-N)^{N-1}$. As is the case with all connected
graphs, the eigenspace associated with the null eigenvalue is the constant vector $\chi_0=\pa{1/\sqrt{N}}[1,...,1]^*$.
Let $\chi_1=\pa{1/\sqrt{2}}[1, -1,0,...,0]$. Then $\inne{\chi_0,\chi_1}=0$ and $L\chi_1=N\chi_1$. Upon solving for the $N-2$ remaining orthonormal
eigenvectors $\chi_l$ for $l=2,...,N-1$, we define the complete graph Fourier transform $\chi_c^*=[\chi_0, \chi_1, \chi_2,...,\chi_{N-1}]^*.$ We then
have $\widehat{\chi_1}=[0, 1, 0,...,0]^*$, and $$\abs{supp(\chi_1)}\abs{supp(\widehat{\chi_1})}=2<N$$ for $N\geq3$; and we see that the support
theorems in
\cite{donoho1989uncertainty} do not hold for graphs. Alternatively, applying Theorem \ref{Koprowski1},  we have, for
$N>2$, that $$\norm{f}^2(N-\sqrt{N})\leq \norm{D_r f}^2+\norm{D_r\widehat{f}}^2\leq \norm{f}^22N.$$ Similarly,
applying Theorem \ref{graph_frame_up}, we have, for $2\leq d\leq N$ and
any $d\times N$ Parseval frame $E$, that
$$2N(d-1)\leq \norm{D_r\chi^*
E^*}^2_{fr}+\norm{D_rE^*}^2_{fr}.$$
    \begin{figure}
\includegraphics[scale=.65]{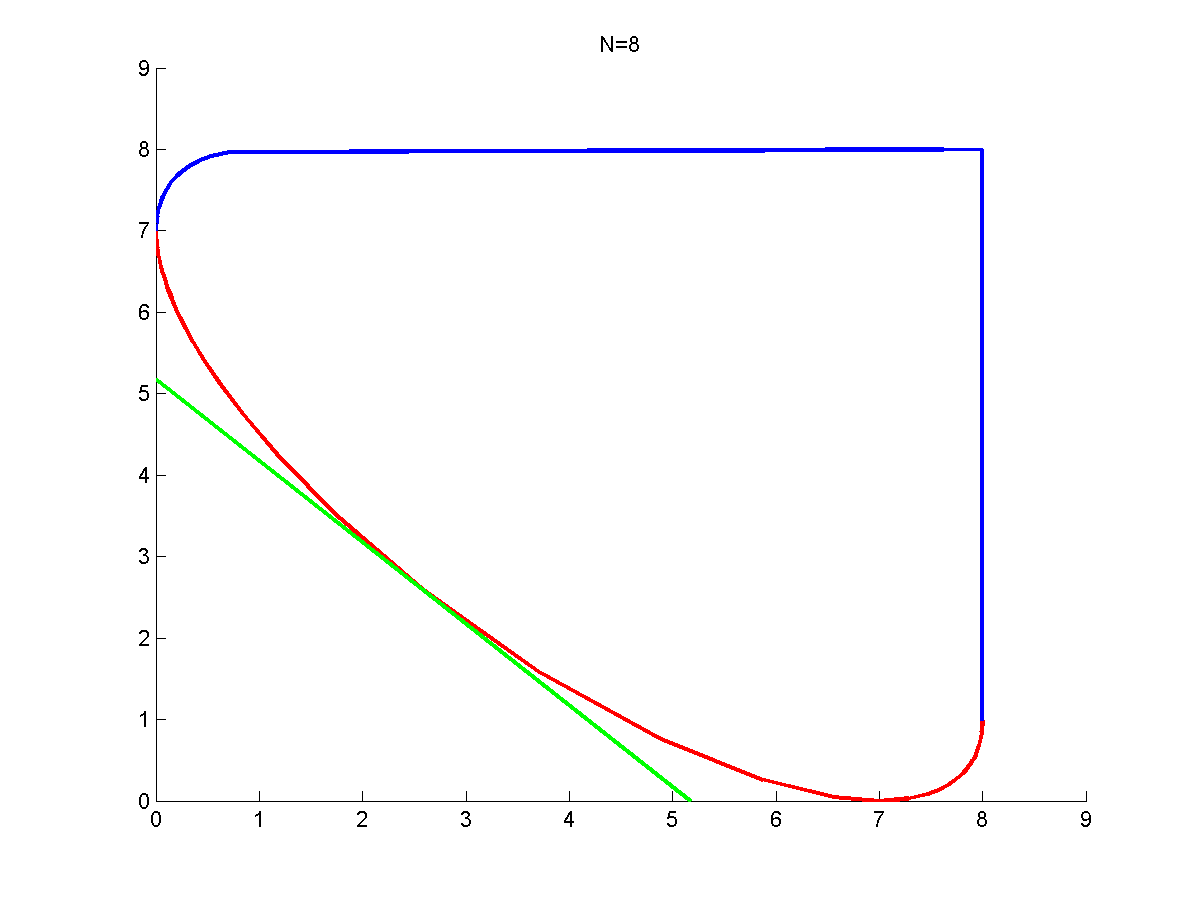}

\includegraphics[scale=.65]{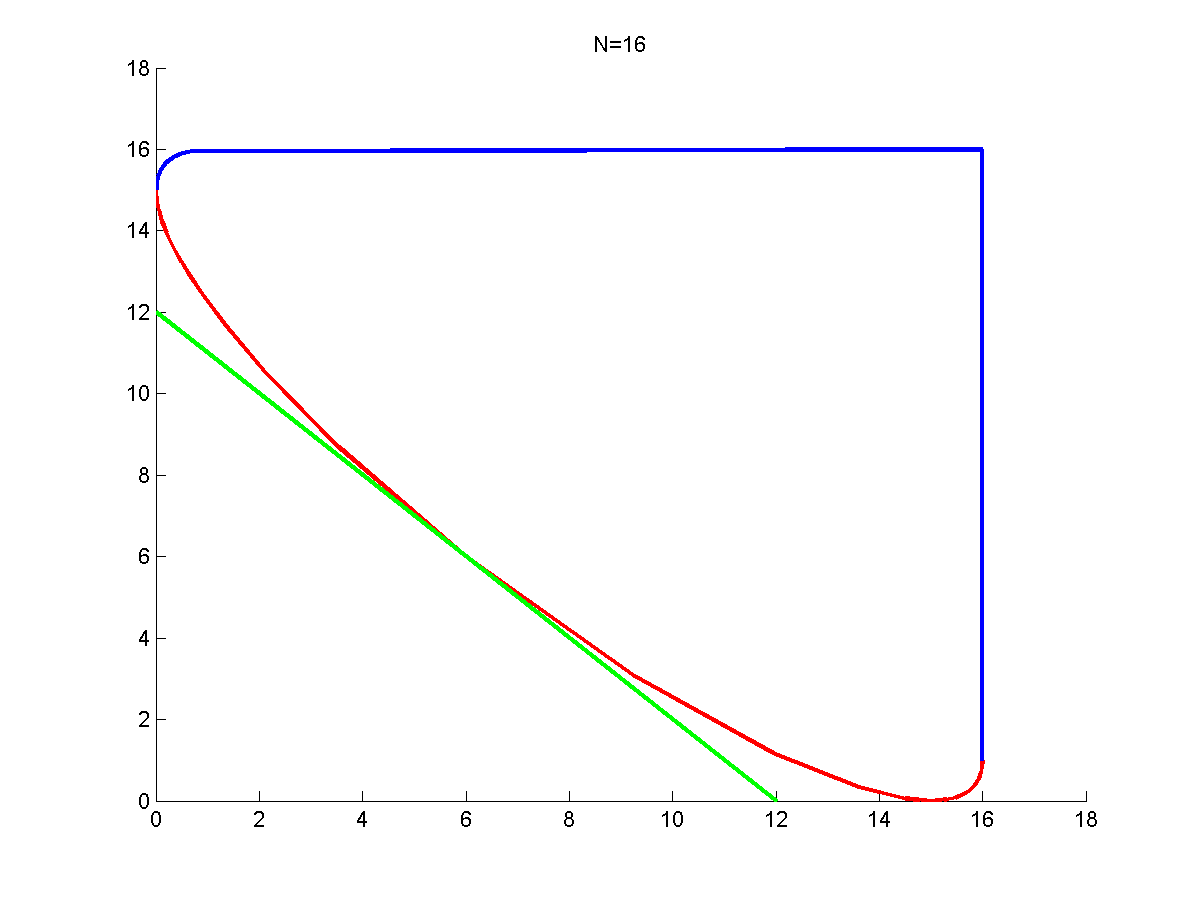}
        \caption{The complete graph differential feasibility regions for 
$N=8$ and $N=16$. The red curve is the differential uncertainty curve,
				the blue is the remaining differential feasibility region boundary, and the
				green line is the line $x+y=N-\sqrt{N}$}\label{figDFR}
\end{figure}

We shall compute the differential feasibility region for the 
complete graph explicitly. We begin our analysis by analyzing the eigenspace of 
$K(\alpha)$.

\begin{prop}\label{prop8}
Let $G$ be the unit weighted complete graph with $N\geq3$ vertices. 
For all $\alpha\neq 0\in\reals$, if $K(\alpha)=L-\alpha \Delta,$
where $L$ is the graph Laplacian for $G$ and $\Delta$ is its diagonalization, then $K(\alpha)$ has an $N-2$ degree
eigenspace associated with the eigenvalue $N(1-\alpha)$.
\end{prop} 

\PROOF $K(\alpha)$ is a block matrix of the form
$$K(\alpha)=\left[
\begin{array}{c|c}
N-1 & -\mathbf{1}_{N-1}^* \\ \hline
-\mathbf{1}_{N-1} & C(\alpha)
\end{array}\right],$$
where $\mathbf{1}_{N-1}$  is the $(N-1)\times 1$ constant function of all ones, and
$C(\alpha)$ is the circulant matrix with $N-1-N\alpha$ on the diagonal
and $-1$ at every other coordinate, i.e., 
$$C(\alpha)=N(1-\alpha)I_{N-1\times N-1}-O_{N-1\times N-1}.$$
Let $W\subset \reals^{N-1}$ be the orthogonal complement of $span(\mathbf{1}_{N-1})$ in $\reals^{N-1}$.
Then for all $w\in W$
we have that $$C(\alpha)w=N(1-\alpha)w-O_{N-1\times N-1}w=N(1-\alpha)w.$$
$W$ has dimension $N-2$ and may be embedded in $l^2(G)$ via the mapping
$w\mapsto \begin{bmatrix} 0, w^*\end{bmatrix}^*.$
We denote this space as $\tilde{W}$. Hence, we have that
$$K(\alpha)\begin{bmatrix} 0\\ w\end{bmatrix}=N(1-\alpha)\begin{bmatrix} 0\\ w\end{bmatrix}$$
and the eigenspace $ES(\alpha)$ associated with $N(1-\alpha)$ has at least dimension $N-2$ 
as it properly contains $\tilde{W}$.
Let $a$ and $b$ denote the remaining two eigenvalues. Let $w_a$ be an eigenvector associated
with $a$ and orthogonal to all  $w\in\tilde{W}$. Then $w_a$ is of the form 
$w_a=c[x 1... 1]^*$ for some real constant $c$. Without loss of generality, we set $c=1$ and we have
$$K(\alpha)w_a=\begin{bmatrix} (N-1)x-(N-1)\\
-x+(1-\alpha N)\\
\vdots \\
-x+(1-\alpha N)\end{bmatrix}=aw_a.$$
Therefore, we must have $a=-x+(1-\alpha N)$. Solving the
quadratic equation resulting from equality in the first coordinate,
i.e., solving the equation, $$x^2-(2-N(\alpha+1))x-(N-1)=0,$$ yields 
$$x=\frac{2-N(\alpha+1)\pm\sqrt{(N(\alpha+1)-2)^2+4(N-1)}}{2}.$$
We conclude that $$a=1-\alpha N-\frac{2-N(\alpha+1)+\sqrt{(N(\alpha+1)-2)^2+4(N-1)}}{2}$$
and $$b=1-\alpha N-\frac{2-N(\alpha+1)-\sqrt{(N(\alpha+1)-2)^2+4(N-1)}}{2}.$$
We conclude that $ES(\alpha)$ has dimension $N-2$ as desired. \hfill $\blacksquare$

From the proof of Proposition \ref{prop8}, we find that the 
minimal eigenvalue of $K(\alpha)$ is 
\begin{align}\lambda_{min}(\alpha)=-\frac{-N(\alpha+1)+\sqrt{(N(\alpha+1)-2)^2+4(N-1)}}{2}-\alpha N,\label{asdff}\end{align}
for all $\alpha\neq 0$. When $\alpha=0$ the minimum eigenvalue is zero, so we may conclude
that equation (\ref{asdff}) holds for all $\alpha\in \reals$.

Let $[x(\alpha), 1,...,1]^*$ with
$$x(\alpha)=\frac{2-N(\alpha+1)+\sqrt{(N(\alpha+1)-2)^2+4(N-1)}}{2}$$
be a vector valued eigenfunction associated with $\lambda_{min}(\alpha)$
for all $\alpha\in \reals$. Upon applying the Rayleigh quotient to 
this vector,
we find that the DUC is the lower boundary of the ellipse
with coordinates 
$$\pa{\frac{N(N-1)}{x(\alpha)^2+(N-1)},\frac{(x(\alpha)-1)^2(N-1)}{x(\alpha)^2+(N-1)}}.$$
The differential feasibility region for $N=8$ and $N=16$ are displayed in
Figure \ref{figDFR}.





\bibliography{sampta_journal_format_2_18}
\bibliographystyle{abbrv}

\vspace{13pt}
\centerline{ACKNOWLEDGEMENT}
\vspace{13pt}

The author gratefully acknowledges the support of the
Wiener Center at the
University of Maryland, College Park and ORAU Maryland.

\end{document}